\documentclass[11pt]{amsart} 
\usepackage{amsmath,amsthm,amsfonts}

\newcommand{\C}{{\mathbb C}}

\newcommand{\Z}{{\mathbb Z}}
\newcommand{\N}{{\mathbb N}}

\newcommand{\OO}{{\mathcal O}}

\newcommand{\BB}{{\mathcal B}}
\newcommand{\DD}{{\mathcal D}}
\newcommand{\EE}{{\mathcal E}}
\newcommand{\CC}{{\mathcal C}}
\newcommand{\AAA}{{\mathcal A}}
\newcommand{\pp}{{\mathfrak p}}
\newcommand{\PP}{{\mathfrak P}}
\newcommand{\qq}{{\mathfrak q}}
\newcommand{\QQ}{{\mathfrak Q}}

\newcommand{\cc}{{\mathfrak c}}
\newcommand{\ttt}{{\mathfrak t}}
\newcommand{\aaa}{{\mathfrak a}}

\newcommand{\ord}{{\operatorname{ord}}}
\newcommand{\Norm}{{\bf{{\operatorname{N}}}}}
\newcommand{\INT}{\operatorname{INT}}
\title[Hilbert's Tenth Problem for Algebraic Function
Fields]{Hilbert's Tenth Problem for Algebraic Function Fields of
  Characteristic $2$} 
\author[Kirsten Eisentr\"ager]{Kirsten Eisentr\"ager}
\address{Department of Mathematics\\ University of California\\
Berkeley, CA 94720,USA}
\email{eisentra@math.berkeley.edu}
\newtheorem{theorem}{Theorem}[section]
\newtheorem{lemma}[theorem]{Lemma}
\newtheorem{cor}[theorem]{Corollary}
\newtheorem{prop}[theorem]{Proposition}
\newtheorem{notation}[theorem]{Notation}
\theoremstyle{definition}
 
\newtheorem{definition}[theorem]{Definition}
\theoremstyle{remark}
\newtheorem{rem}{Remark$\!\!$}        
\date{November 3, 2002} 
\begin{document}
\begin{abstract}
  {\small Let $K$ be an algebraic function field of characteristic $2$
    with constant field $C_K$. Let $C$ be the algebraic closure of a
    finite field in $K$. Assume that $C$ has an extension of degree
    $2$. Assume that there are elements $u,x$ of $K$ with $u$
    transcendental over $C_K$ and $x$ algebraic over $C(u)$ and such
    that $K=C_K(u,x)$. Then Hilbert's Tenth Problem over $K$ is
    undecidable. Together with Shlapentokh's result for odd
    characteristic this implies that Hilbert's Tenth Problem for any
    such field $K$ of finite characteristic is undecidable. In
    particular, Hilbert's Tenth Problem for any algebraic function
    field with finite constant field is undecidable.}
  \end{abstract}
\maketitle
\section{Introduction}
Hilbert's Tenth Problem in its original form can be stated in the
following form: Is there a uniform algorithm that determines, given a
polynomial equation with integer coefficients, whether the equation
has an integer solution or not? In \cite{Mat70} Matijasevich proved
that the answer to this question is no, {\it i.e.}\ that Hilbert's
Tenth Problem is undecidable. Since then various analogues of this
problem have been studied by asking the same question as above for
polynomial equations with coefficients and solutions over some other
commutative ring $R$.  Perhaps the most important unsolved question in
this area is Hilbert's Tenth Problem over the field of rational
numbers. There are also many results that prove undecidability: It was
proved in \cite{Den80} and \cite{DenLi78} that Hilbert's Tenth Problem
is undecidable for various rings of algebraic integers, and
\cite{Den78} proves the undecidability of the problem for rational
functions over formally real fields. In \cite{KR92} Kim and Roush
proved that Hilbert's Tenth Problem over $\C(t_1,t_2)$ is undecidable.
Diophantine undecidability has also been proved for some rational
function fields of characteristic~ $p$: Pheidas \cite{Phei91} has
shown that Hilbert's Tenth Problem is undecidable for rational
function fields over finite fields of characteristic greater than~$2$
and Videla \cite{Vi94} has proved the analogous result for
characteristic $2$. Kim and Roush \cite{KR} proved undecidability for
rational function fields of characteristic $p>2$ whose constant fields
do not contain the algebraic closure of a finite field. In
\cite{Sh2000} Shlapentokh proved that the problem for algebraic
function fields over possibly infinite constant fields of
characteristic $p>2$ is undecidable. This paper will solve the
analogous problem over function fields of characteristic $2$, so
Hilbert's Tenth Problem for any such field of finite characteristic is
undecidable.  We will first describe the general approach that is used
to prove the undecidability of Hilbert's Tenth Problem for any
function field of positive characteristic. The approach is based on an
idea that was first introduced by Denef in \cite{Den79} and further
developed by Pheidas in \cite{Phei91} and Shlapentokh in \cite{Sh96}
and \cite{Sh2000}.

Before we can describe the idea in detail we need to define what an
algebraic function field is:

\begin{definition}
  A field extension $K/C_K$ is said to be an {\em algebraic function field
  (of one variable)} if these conditions hold:
\begin{enumerate}
\item the transcendence degree of $K/C_K$ is $1$.
\item K is finitely generated over $C_K$; and
\item $C_K$ is algebraically closed in $K$.
\end{enumerate}
In this case there exists $t \in K$, transcendental over $C_K$, such
that the degree of the field extension $[K:C_K(t)]$ is finite. The
field $C_K$ is called the constant field of $K$.
\end{definition}
We also need to define two notions that we will use below:
\begin{definition}
  1. If $R$ is a commutative ring, a {\em diophantine equation over
    $R$} is an equation $P(x_1,\dots,x_n)=0$ where $P$ is a polynomial
  in the variables $x_1, \dots, x_n$ with coefficients in
  $R$.\\
  2. A subset $S$ of $R^k$ is diophantine if there is a polynomial
  $P(x_1,\dots,x_k,$ $y_1,\dots,y_m) \in
  R[x_1,\dots,x_k,y_1,\dots,y_m]$ such that \[S=\{(x_1, \dots, x_k)
  \in R^k: \exists y_1, \dots,y_m \in R,\;
  (P(x_1,\dots,x_k,y_1,\dots,y_m)=0)\}.\]
\end{definition}
When $R$ is not a finitely generated algebra over $\Z$, we restrict our
attention to diophantine equations whose coefficients are in a
finitely generated algebra over $\Z$. In particular, if $R$ is a ring
of polynomials or a field of rational functions in an indeterminate
$t$, we only consider diophantine equations whose coefficients lie in
the natural image of $\Z[t]$ in $R$.
%%%%%%%%%%%%%%%%%%%%%%%%%%%%%%%%%%%%%%%%%%%%%%
%%%%%%%%%%%%%%%%%%%%%%%%%%%%%%%%%%%%%%%%%%%%%%%
%%%%%%%%%%%%%%%%%%%%%%%%%%%%%%%%%%%%%%%%%%%%%%%
\subsection{Idea of Proof.} 
Let $\N$ be the set of natural numbers $\{0,1,2,\dots\}.$ The general
idea of the proof is to reduce a certain decision problem over the
natural numbers which we know to be undecidable to Hilbert's Tenth
Problem over $K$. The undecidable structure that we will use is the
diophantine theory of the natural numbers with addition and a
predicate $|_p$ defined by $n|_pm $ if and only if $\exists s \in \N
(m = p^s n)$. In \cite{Phei87} Pheidas showed that this structure has
an undecidable diophantine theory, {\it i.e.}\ there is no uniform
algorithm that, given a system of equations over the natural numbers
with addition and $|_p$, determines whether this system has a solution
or not. To reduce this problem to Hilbert's Tenth Problem over $K$ we
first let $G$ be a subfield of $K$ containing an element $t$
transcendental over $C_K$. The field $G$ will be defined in
Lemma~\ref{first}. Also fix a prime $\pp$ of $K$ which lies above a
nontrivial prime of $G$. We can choose $t$ and $\pp$ such that
$\ord_{\pp}t=1$. Both $t$ and $\pp$ will be defined at the end of
section~\ref{setup}. Let $\OO_{K,\pp}:=\{x \in K: \ord_{\pp}x \geq
0\}$, and let $\OO_{G,\pp}:= G \cap \OO_{K,\pp}$.  Now let $\INT(\pp)$
be any subset of $K$ such that $\OO_{G,\pp}\subseteq
\INT(\pp)\subseteq \OO_{K,p}.$ We define a map $f$ from the integers
to subsets of $K$ by associating to an integer $n$ the subset
$f(n):=\{x \in \INT(\pp): \ord_{\pp}x =n\}$.  Then $n_3=n_1+n_2$ ($n_i
\in \N$) is equivalent to the existence of $z_i \in f(n_i)$ such that
$z_3=z_1\cdot z_2.$ This follows from the fact that $\ord_{\pp}z_1 +
\ord_{\pp}z_2=\ord_{\pp}(z_1 \cdot z_2)$ and that $t^{n_i} \in
f(n_i)$.  We also have for natural numbers $n,m$
\begin{eqnarray*}
n|_p m &\iff& \exists s \in \N\;\; m = p^s n \\&\iff&
\exists x \in f(n)\; \exists y \in f(m)\; \exists s \in \N \;\; 
(\ord_{\pp}y=p^s \ord_{\pp}x).
\end{eqnarray*}
This equivalence can be seen easily, because we can let $x:=t^n$ and $y:=t^m$.
But the last formula is equivalent to
\[
\exists x \in f(n)\; \exists y \in f(m)\; \exists w \in K \;\exists s
\in \N\; w=x^{p^s} \mbox{ and } \{w/y,y/w\} \subset \INT(\pp) .\]
Here $w/y\in \INT(\pp)$ and $y/w \in \INT(\pp)$ just means that $y$
and $w$ have the same order at $\pp$.

If we have diophantine definitions for $p(K):=\{(x,w) \in K^2: \exists
s \in \N, w=x^{p^s}\}$ and $\INT(\pp)$, then the above argument shows
that for every system of equations with addition and $|_p$ we can
construct a system of polynomial equations over $K$ which will have
solutions in $K$ if and only if the original system of equations over
$\N$ has solutions in $\N$. But the diophantine theory of $\N$ with
$+$ and $|_p$ is undecidable; hence Hilbert's Tenth Problem over $K$
is undecidable.

So the strategy for the proof will be to prove that $p(K)$ is
diophantine and that there exists some set $\INT(\pp)$ as above which
is diophantine for the class of fields $K$ that we are
considering. This can be summarized as

\begin{theorem} \label{THEOREM}
  Let $K$ be an algebraic function field of characteristic $2$ with
  constant field $C_K$. Let $C$ be the algebraic closure of a finite
  field in $K$. Assume that $C$ has an extension of degree $2$. Assume
  that there are elements $u,x$ of $K$ with $u$ transcendental over
  $C_K$ and $x$ algebraic over $C(u)$ and such that $K=C_K(u,x)$. Then
  $p(K)$ is diophantine. Also there exists a subfield $G$ of $K$ as
  above with $C(t) \subseteq G$ for an element $t$ transcendental over
  $C_K$.  There exists a prime $\pp$ of $K$ satisfying the conditions
  above such that $\INT(\pp)$ is diophantine for some set $\INT(\pp)$
  with $\OO_{G,\pp} \subseteq \INT(\pp) \subseteq \OO_{K,\pp}$. So
  Hilbert's Tenth Problem over $K$ is undecidable.
\end{theorem}
In \cite{Sh2000} Shlapentokh proves that for such $K$ in any
characteristic $p>0$ there exists some set $\INT(\pp)$ as above which
is diophantine. She also proves that $p(K)$ is diophantine when the
characteristic of $K$ is greater than $2$, but her main lemmas are not
valid in characteristic $2$. So in order to prove undecidability in
characteristic $2$, the last open case, we need to prove that $p(K)$
is diophantine when the characteristic of $K$ is $2$. The rest of the
paper is devoted to proving this. The outline of the proof follows
Shlapentokh's proof for odd characteristic. Before we can prove this
we first need to prove some properties of $K$ and then set up some
notation. The next section will do that. In section~\ref{p-power} we
will prove that the set $p(K)$ is diophantine in characteristic $2$.
\section{Setup and Notation}\label{setup}
Let $\N$ be the set of natural numbers $\{0,1,2,\dots\}.$ Let
$K,C_K,C,u$ and $x$ be as in Theorem~\ref{THEOREM}.  We will use the
following
\begin{notation}
Let $F$ be a field, and $k \in \N$. We denote by $F^k$ the set
$ F^{k}:= \{a^k:a \in F\}.$
\end{notation}
We will now prove some properties of $K$ that we will need later on.
We may assume that $u$ is not a square in $K$, because if $u
= u_1^{2}$ with $u_1 \in K$ and $s \in \N$, we can replace $u$ by
$u_1$. Then $K=C_K(u,x)=C_K(u_1,x)$. Since the extension $K/C_K(u)$
can be generated by a single element, $u \in K^{2^s}$ only if $s \leq
[K:C_K(u)]$, so replacing $u$ by its square root terminates after a
finite number of steps.
\\

We have the following
\begin{lemma}\label{first}  Let $K,C_K,C,u,x$ be as above.
  Let $G$ be the algebraic closure of $C(u)$ inside $K$.  Then
  $G=C(u,x)$.
\end{lemma}
\begin{proof}
  First note that $C(u)$ is algebraically closed in $C_K(u)$, because
  $C$ is algebraically closed in $C_K$ (\cite{Deu73}, p.~117).
  Let $m:=[K:C_K(u)]$. If $m=1$, {\it i.e.}\ $x \in C_K(u)$, then the
  statement is true since $C(u)$ is algebraically closed in $C_K(u)$.
  So assume $x \notin C_K(u)$. Let $\alpha \in G$, $\alpha \notin
  C(u)$. Then by \cite{La} Lemma~4.10, p.~366,
  \begin{equation} \label{inequ}
[C(u,\alpha):C(u)]=[C_K(u,\alpha):C_K(u)] \leq [K:C_K(u)]=m.
\end{equation}
In particular, $[C(u,x):C(u)] = m$. Now assume by contradiction that
there exists a $\beta \in G$, $\beta \notin C(u,x)$. Let $G_1:=
C(u,x,\beta)$. Then $[G_1:C(u)] > m$. Also $G_1$ is an algebraic
function field with constant field $C$, and $C$ is perfect. Then by
\cite{Ma84}, p.~94 the extension $G_1/C(u)$ is finite and separable,
since $u$ is not a square in $G_1$. Hence there exists a primitive
element $\gamma \in G_1$ with $C(u,\gamma)=G_1$. But then
$[C(u,\gamma):C(u)]=[G_1:C(u)]>m$, contradicting (\ref{inequ}).
\end{proof}
%%%%%%%%%%%%%%%%%%%%%%%%%%%%%%%%%%%%%%%%%%%%%%%%%%%%%%%%%%%%
%%%%%%%%%%%%%%%%%%%%%%%%%%%%%%%%%%%%%%%%%%%%%%%%%%%%%%%%%%%%
\begin{definition}
Let $K$ be an algebraic function field with constant field $C_K$. A
{\em constant field extension of $K$} is an algebraic function field $L$
with constant field $C_L$ such that $L \supseteq K$, $C_L \cap K =
C_K$ and $L$ is the composite extension of $K$ and $C_L$, $L=C_LK$.
\end{definition}
%%%%%%%%%%%%%%%%%%%%%%%%%%%%%%%%%%%%%%%%%%%%%%%%%%%%%%%%%%%%%%%%
\begin{prop}\label{existence}
  Let $G$ be as in Lemma~\ref{first}.  Fix a positive integer $k$.
  For any sufficiently large positive integer $h$ a finite constant
  extension of $G$ contains a nonconstant element $t$ and a set of
  constants $V$ of cardinality $k+2^h$ such that $0\in V$, $1 \notin
  V$. Also we can choose $t$ and $V$ such that for all $c \in V$ the
  divisor of $t+c$ is of the form $\pp_c/\qq$, where the $\pp_c$'s and
  $\qq$ are prime divisors of degree $2^h$.
\end{prop}
\begin{proof}
  This is Theorem~6.11 of \cite{Sh2000} if $C$ is infinite. The proof
  of the existence of $t$ and $V$ with the desired properties
  in Theorem~6.11 does not use that $C$ is infinite; it only
  requires passing to a finite extension of $C$.
\end{proof}
\begin{rem} In Proposition~\ref{existence} we can choose $V$ with the
  property that for all $s \in \N$ for all $c,c'\in V$ $c^{p^s}\neq
  c'$ if $c\neq c'$.
\end{rem}
From now on we will assume that an element $t$ and a set $V$ of
constants with the desired properties as in
Proposition~\ref{existence} already exist in $G$. (Otherwise rename
the constant extension $G$ again and work with it instead.)  Enlarging
the field of constants by a finite extension is okay as far as the
undecidability of Hilbert's Tenth Problem is concerned. Also let $\pp
:= \pp_0$, so that the divisor of $t$ is of the form $\pp/\qq$.
\begin{prop}\label{2h}
  Let $G,C,t$ be as above. Then $[G:C(t)]$ is separable, and $2^h =
  n=[G:C(t)]$. 
\end{prop}
\begin{proof}
  Since the divisor of $t$ is of the form $\pp/\qq$, $t$ is not a
  square in $G$. Also $C$ is perfect. Hence $G/C(t)$ is separable by
  \cite{Ma84}, p.~22. Also by \cite{FriJa86}, p.~13, $[G:C(t)]= \deg
  \pp = \deg \qq=2^h$.
\end{proof}
Now we can prove that $K$ is separably generated:
\begin{cor}
  Let $K$ be an algebraic function field with constant field $C_K$.
  Let $C$ be the algebraic closure of a finite field in $K$. Assume
  that $C$ has an extension of degree $2$. Assume that there exist
  $x,u$ as above. Let $G,t$ be as above. Then $K/C_K(t)$ is separable.
\end{cor}
\begin{proof}
  By Proposition \ref{2h} $G/C(t)$ is separable. The field $K$ is the
  compositum of $C_K(t)$ and $G$ over $C(t)$, hence $K/C_K(t)$ is also
  separable.
\end{proof}
Now we can use Lemma~6.13 of \cite{Sh2000} to see how the $\pp_c$'s
and $\qq$ behave in the extension $K$.
\begin{lemma}[Lemma~6.13 of \cite{Sh2000}]\label{Sh}
   Let H be an algebraic function field over a field of constants
   $C_H$.  Let $K$ be a constant field extension of $H$. Let $C_K$ be
   the constant field of $K$, and assume $H$ is algebraically closed in
   $K$. Let $t \in H-C_H$ be such that $H/C_H(t)$ is separable. Let
   $\aaa$ be a prime of $C_H(t)$ remaining prime in the extension $H$
   and such that its residue field is separable over $C_H$. Then
   $\aaa$ will have just one prime factor in $K$.
\end{lemma}
This lemma easily implies the following corollary:
\begin{cor}
  Let $\{\pp_c: c\in V \}$ and $\qq$ be as in Proposition
  \ref{existence}. Then the $\pp_c$'s and $\qq$ remain prime in $K$.
\end{cor}
\begin{proof}
  Lemma \ref{Sh} applies, since $K/C_K$ is a constant
  field extension of $G/C$: by construction $C$ is algebraically
  closed in $C_K$, and also $C_KG=K$. The only thing we need to check
  is that $G \cap C_K=C$. Assume $\alpha \in G-C$. Then $\alpha$ is
  transcendental over $C$ and also over $C_K$. Hence $\alpha \notin G
  \cap C_K$. Thus we can apply the lemma to the primes, $t,1/t$ and
  $t+c$ of $C(t)$. Since $C$ is perfect, the residue extensions of
  the primes will be separable.
\end{proof}
Since the $\pp_c$'s and $\qq$ remain prime in $K$ we will just denote
them by the same letters again when considering them as primes of $K$,
and we will let $\pp:=\pp_0$.  Now we can fix some notation that we
will use for the rest of the paper:
\begin{itemize}
\item $K$ will denote an algebraic function field over a field of
  constants $C_K$ of characteristic $p =2$.
\item $C$ will denote the algebraic closure of a finite field inside
  $C_K$.
\item $t$ will denote a nonconstant element of $K-C_K$ such that the
  divisor of $t$ is of the form $\pp/\qq$, where $\pp$, $\qq$ are
  primes of degree $2^h$ for some natural number $h$.
  Furthermore, $K/C_K(t)$ is separable, and $2^h =n = [K:C_K(t)]$.
\item $\tilde C_K$ will denote the algebraic closure of $C_K$, and
  $\tilde K := \tilde C_K K$.
\item $r$ will denote the number of primes of $\tilde K$ ramifying in
  the extension $\tilde K/\tilde C_K(t)$.
\item $V$ will denote a subset of $C$, containing $n+2r+6$
  elements, such that $0 \in V$, $1 \notin V$, and such for all $ c\in
  V$ the divisor of $t+c$ is of the form $\pp_c/q$, where $\pp_c$ is a
  prime divisor of $K$. Also pick $V$ such that for any $s \in \N$,
  $c,c' \in V$, we have $c^{p^s} \neq c'$, if $c \neq c'$.
\item For all $c \in V$, $\PP_c$ will denote the prime of $C_K(t)$
  lying below $\pp_c$, while $\QQ$ will denote the prime of $C_K(t)$
  lying below $\qq$. Also let $\PP:=\PP_0$. For all $c \in V$, $\PP_c$
  and $\QQ$ do not split in the extension $K/C_K(t)$.
\item For every $c \in V$, $V_c$ will denote the set $V_c:= \{c^{p^j}:
  j \in \N \}.$ Since every $c \in V$ is algebraic over a finite
  field, $V_c$ is a finite set for all $c \in V.$

\end{itemize}
To obtain $t$ and $V$ with the desired properties, we
have to assume that $C$ is sufficiently large, but this is not a
restriction because we can enlarge the field of constants and by
Proposition \ref{existence} a finite extension is enough. Let $L$ be
this finite extension. If Hilbert's Tenth Problem over $L$ is
undecidable, then Hilbert's Tenth Problem over $K$ is also
undecidable.  So in the following we will assume that $L=K$ to
simplify notation.

\section{$p$-th Power Equations}\label{p-power}
Using the notation that we set up in the last section will now prove
that the set $p(K)=\{(x,y) \in K^2 :\exists s \in \N, y=x^{2^s}\}$ is
diophantine which is Theorem~\ref{THEEND} below.  The main ingredient
for proving this is the next theorem. It gives an equivalent
definition of what it means for $(x,y)$ to be in $p(K)$. Eventually we
want to find polynomial equations describing these relations, so the
goal afterwards will be to rewrite the equations below as polynomial
equations.
\begin{theorem}\label{mainpart}
  Given $x, y \in K$, let $u:= \frac{x^2 + t^2 + t}{x^2 + t}$ and
  $\tilde u := \frac{x^2 + t^ {-2} + t^{-1}}{x^2 + t^{-1}}$.
Let $v := \frac{y^2 + t^{2^{s+1}}+ t^{2^s}}{y^2 + t^{2^s}}$ and
$\tilde v := \frac{y^2 + t^{-2^{s+1}}+ t^{-2^s}}{y^2 +
  t^{-2^s}}$ for some $s \in \N.$

Then $y=x^{2^s}$ if and only if
\begin{gather}
 \label{E1:int}  \exists r \in \N \;v = u^{2^r}\\ 
\label {E2:int} \exists j \in \N \;\tilde v = \tilde u^ {2^j}. 
\end{gather}
\end{theorem}
\begin{proof}
  Suppose $y=x^{2^s}$. Let $r=j=s$. Then (\ref{E1:int}) and
  (\ref{E2:int}) are satisfied. This completes one direction of the
  proof.

On the other hand, suppose that $r$ and $j$ as in the statement of the
theorem exist. Then 
\[
v = \left( \frac{x^2 + t^2 + t}{x^2 +t} \right) ^{2^r} =
\frac{x^{2^{r+1}}+ t^{2^{r+1}}+ t^{2^r}}{x^{2^{r+1}}+t^{2^r}} =
\frac{y^2 + t^{2^{s+1}}+t^{2^s}}{y^2 + t^{2^s}}.
\]
So
\[
(x^{2^{r+1}}+ t^{2^{r+1}}+ t^{2^r})(y^2 + t^{2^s})=(x^{2^{r+1}}+
t^{2^r})(y^2 + t^{2^{s+1}}+t^{2^s}),
\]
{\it i.e.}\
\[
t^{2^{r+1}}y^2 + t^{2^{r+1}+ 2^s} =
x^{2^{r+1}}t^{2^{s+1}}+t^{2^r+2^{s+1}}.
\]
Thus
\begin{gather}\label{G:int}
y^2 = (x^{2^{r+1}}t^{2^{s+1}}+t^{2^r + 2^{s+1}}+t^{2^{r+1}+2^s})\cdot
t^{-2^{r+1}}.
\end{gather}
Hence if we can show that $r=s$, then $y^2 = x^{2^{s+1}}$,
so  $y=x^{2^s}$, since the characteristic of $K$ is 2. So our goal is to show
that $r=s$.

Similarly to the calculations above we get
\[
\tilde v = \left( \frac{x^2 + t^{-2} + t^{-1}}{x^2 +t^{-1}} \right)
^{2^j} = \frac{x^{2^{j+1}}+ t^{-2^{j+1}}+
  t^{-2^j}}{x^{2^{j+1}}+t^{-2^j}} = \frac{y^2 +
  t^{-2^{s+1}}+t^{-2^s}}{y^2 + t^{-2^s}},
\]
and we get

\begin{gather} \label{H:int}
y^2 = (x^{2^{j+1}}t^{-2^{s+1}}+t^{-2^j - 2^{s+1}}+t^{-2^{j+1}-2^s})\cdot
t^{2^{j+1}}.
\end{gather}

By (\ref{G:int}) \begin{gather}\label{I:int}
y= (x^{2^r}t^{2^s} +
t^{2^{r-1}+2^s}+t^{2^r+2^{s-1}})\cdot t^{-2^r} \end{gather}
 (unless $r$ or $s$ are $< 1$), and by (\ref{H:int})
\begin{gather} \label{J:int}
y=(x^{2^j}t^{-2^s}+t^{-2^s-2^{j-1}}+t^{-2^j-2^{s-1}})\cdot t^{2^j}
\end{gather} (unless $j$ or $s <1$). Eliminating $y$ from (\ref{I:int})
and (\ref{J:int}), we get
\begin{gather}\label{K:int}
  (t^{2^s-2^r}x^{2^r})+(t^{2^j-2^s}x^{2^j})=
  t^{2^s-2^{r-1}}+t^{2^{s-1}} + t^{2^{j-1}-2^s}+t^{-2^{s-1}}.
\end{gather}

Now assume that $y$ is a square, say $y=z^2$ (and $s,j,r >0$). Then
\[
v= \frac{z^4 + t^{2^{s+1}}+t^{2^s}}{z^4 + t ^{2^s}} = \left( \frac
  {z^2 + t^{2^s}+t^{2^{s-1}}}{z^2 + t^{2^{s-1}}} \right) ^2 = (v')^2.
\] 
Hence
\[
v = (v')^2=u^{2^r}\text{ , so } u^{2^{r-1}}= \left(
  \frac{z^2+t^{2^s}+t^{2^{s-1}}}{z^2 + t^{2^{s-1}}} \right) = v'.
\]
Similarly $\tilde v = (\tilde {v}')^2$ and $\tilde u^{2^{j-1}}= \tilde
{v}'$, so in the new formulae $s,r$ and $j$ are replaced by $s-1,r-1$
and $j-1$, respectively, and we're done if we can show that $z =
x^{2^{s-1}}$. Hence we can reduce the problem to the case where either
(a) $s=0$ or $r=0$ or $j=0$, or (b) $y$ is not a square.

Case (a): $s=0$: If $s=0$, then $v=\frac{y^2+t^2+t}{y^2+t}$, and $v$ is
not a square since $\frac{dv}{dt}= \frac{y^2 + t^2 + t+ y^2 +t}{(y^2
  +t)^2} = \frac{t^2}{(y^2 +t)^2}\neq 0$. So if $s=0$, then $v=
\frac{y^2 + t^2 +t}{y^2 + t} = u^{2^r}$. Since $v$ is not a square,
this implies $r=0$. Hence $r=s=0$ and we're done.

If $r=0$, then $v=u$. By the same argument as above $u$ is not a
square. Now if $s>0$, then $v$ is a square and hence $u$ is a square,
contradiction. Hence $r=s=0$, and we're done. The case $j=0$ follows
from symmetry.

Case (b): By case (a) we may assume $r>0, s>0$ and $j>0$ and by
contradiction let's assume that $r \neq s$. If we look at equations
(\ref{I:int}) and (\ref{J:int}), we see that $y$ is a square unless
(i) $s=1$ or (ii) both $r=j=1$.

(i) Suppose $s=1$.  Since we're done if $r=s$ we may assume that $r
\geq~2, j\geq 2$. From (\ref{K:int}) we obtain
\[
t^{2-2^r}x^{2^r}+ t^{2^j-2}x^{2^j}=t^{2-2^{r-1}}+t^{2^{j-1}-2}+t +
\frac{1}{t}
\]
or
\[
(t^{1-2^{r-1}}x^{2^{r-1}}+
t^{2^{j-1}-1}x^{2^{j-1}}+t^{1-2^{r-2}}+t^{2^{j-2}-1})^2 = t+ \frac{1}{t}
\]
Since $j \geq 2$ and $r \geq 2$ the left side is a square. The right
side is not, contradiction.

(ii) Suppose $r=j=1$. Again since we're done if $r=s$ we may assume
$s>1$.  By (\ref{K:int}) we have
\[
x^2(t^{2^s-2}+t^{2-2^s})=t^{1-2^s}+t^{2^s-1} + \frac{1}{t^{2^{s-1}}} +
t^{2^{s-1}}.
\]
Let $\pp$ be the simple zero of $t$. Since $1-2^s < -2^{s-1}$ $(s \geq
2)$, the right side has a pole of odd order at $\pp$, while the left
side is a square, so it only has poles of even order. This proves the
theorem.
\end{proof}
So the goal for the rest of this section is to show that the relations
we used in the statement of the theorem are diophantine. To do that it
will clearly be enough to show that the following four sets are
diophantine:
\[ S:=\{\,t^{2^s}: s \in \N \,\},\;\;\;S':=\{\,{(t^{-1})}^{2^s}:
s \in \N \,\},\]
\[
T:= \left\{(x,w)\in K^2 : \exists s \in \N, w =\left(
  \frac{x^2+t^2+t}{x^2+t} \right)^{2^s}\right\},
\]
and
\[ T':= \left\{(x,w)\in K^2 : \exists s \in \N, w =\left(
  \frac{x^2+(t^{-1})^2+t^{-1}}{x^2+t^{-1}} \right)^{2^s}\right\}\]

It is enough to prove that $S$ and $T$ are diophantine, because we can
replace $t$ by $t^{-1}$ and replace $V$ by $W:= \{\frac{1}{c}:c\in
(V-\{0\})\}\cup\{0\}$ in section~\ref{setup}. Then we can also replace
$t$ by $t^{-1}$ and $V$ by $W$ in the whole proof to obtain
diophantine definitions for $S'$ and $T'$.

Lemma~\ref{tpower} and Corollary~\ref{diophantine1} below will show
that $S$ is diophantine, and Corollary~\ref{diophantine2} will show
that $T$ is diophantine.\\
\subsection{\bf \large {The set $S=\{\,t^{2^s}: s \in \N \,\}$
    is diophantine.\\}}

To prove that $S$ is diophantine, we first need a definition and a
lemma:
\begin{definition}
  Let $w \in K$. The {\em height} of $w$ is the degree of the zero
  divisor of $w$.
\end{definition}
\begin{rem}Equivalently, we could have defined the height of $w \in K$
  to be the degree of the pole divisor of $w$.\end{rem}
%%%%%%%%%%%%%%%%%%%%%%%%%%%%%%%%%%%%%%%%%%%%%%%%%%%%%%%%%%%%%
%%%%%%%%%%%%%%%%%%%%%%%%%%%%%%%%%%%%%%%%%%%%%%%%%%%%%%%%%%%%%
%
%
\begin{lemma}\label{frac}
  Let $w \in K$, let $a,b \in C$. Then all the zeros of $\frac{w
    +a}{w+b}$ are zeros of $w+a$ and all the poles of
  $\frac{w+a}{w+b}$ are zeros of $w+b$. Furthermore, the height of
  $\frac{w+a}{w+b}$ is equal to the height of $w$.
\end{lemma}
\begin{proof}
This is Lemma~2.4 in \cite{Sh2000}.
\end{proof}

%%%%%%%%%%%%%%%%%%%%%%%%%%%%%%%%%%%%%%%%%%%%%%%%%%%%%%%%%%%%%%%%%
%%%%%%%%%%%%%%%%%%%%%%%%%%%%%%%%%%%%%%%%%%%%%%%%%%%%%%%%%%%%%%%%%

\begin{lemma}\label{basecase}
  Let $u,v,z \in \tilde K:= \tilde C_K K$, assume that $z\notin \tilde
  C_K$, and let $y \in \tilde C_K(z)$. Assume that $y,z$ do not have
  zeros or poles at any valuation of $\tilde K$ ramifying in the
  extension $\tilde K/ \tilde C_K(z)$ and that $\tilde K/ \tilde
  C_K(z)$ is separable.  Moreover, assume
\begin{gather}
 \label{eqn} y+z = u^4 + u\\
  \frac{1}{y} + \frac{1}{z} =
  v^4 + v
\end{gather}
Then $y = z ^{4^k}$ for some $k\geq 0$.
\end{lemma}
\begin{proof}
  Recall that for a field $F$ and a natural number $k$, $F^k$
  denotes the set $F^k=\{a^k:a \in F\}.$ 
  In $\tilde C_K(z)$ the zeros and poles of $z$ are
  simple. Assuming that $z$ satisfies the conditions of
  Lemma~\ref{basecase} thus amounts to assuming that all zeros and
  poles of $z$ are simple in $\tilde K$.
  
  Equation~(\ref{eqn}) and the fact that $z$ has simple poles imply
  that $y \notin \tilde C_K$, so $y \in (\tilde K)^{4^s}$ only if $s
  \leq [\tilde K :\tilde C_K(z)]$. If $y=w^4$ with $w \in \tilde
  C_K(z)$, then $w+z = (u+w)^4 +(u+w)$ and $1/w +1/z=
  (v+1/w)^4+(v+1/w)$. So if we can prove that $w=z^{4^s}$ for some $s
  \in \N$, then $y=w^4 = z^{4^{s+1}}$. Hence we may assume without
  loss of generality that $y \notin (\tilde C_K(z))^4$.
 
  Let $\frac{\AAA}{\BB}$ be the divisor of $z$ in $\tilde K$,
  where $\AAA$ and $\BB$ are relatively prime effective divisors.  By
  assumption, all the prime factors of $\AAA$ and $\BB$ are distinct.
  Also all the poles of $u^4 +u$ and $v^4+v$ have orders
  divisible by $4$.
  
  {\it Claim:} The divisor of $y$ is of the form $\EE^4 \DD$ where all
  the prime factors of $\DD$ come from $\AAA$ or $\BB$.  Also the
  factors of $\AAA$ that appear in $\DD$, will appear to the first
  power in $\DD$ and the factors of $\BB$ that appear in $\DD$ occur
  to the power $-1$.\\ {\it Proof of Claim:} Let $\ttt$ be a prime
  which is not a factor of $\AAA$ or $\BB$. Without loss of generality
  assume $\ttt$ is a pole of $y$. Then, since $\ord_{\ttt}z =0$, we
  have
\[
0 > \ord_{\ttt}y = \ord_{\ttt}(z + y) = \ord_{\ttt}(u^4 + u) \equiv 0
\mod 4.
\]
Now let $\ttt$ be a factor of $\AAA$ or $\BB$. Again without loss of
generality assume $\ttt$ is a pole of $y$. If $\ttt$ is a factor of
$\AAA$, then $\ord_{\ttt}y = \ord_{\ttt}(y+z) = \ord_{\ttt}(u^4 + u)$.
Hence $\ttt$ is a pole of $u$, so $\ord_{\ttt}y \equiv 0 \mod 4$. If,
however, $\ttt$ is a factor of $\BB$, there are two possibilities:
either $\ord_{\ttt}y = \ord_{\ttt}z=-1$ or again
$\ord_{\ttt}y=\ord_{\ttt}(u^4+u) \equiv 0 \mod 4$. This proves the
claim.

On the other hand, $\AAA$ and $\BB$ considered as divisors over
$\tilde C_K(z)$ are prime divisors, and since $y \in \tilde C_K(z)$,
we can deduce that the divisor of $y$ is of the form $\EE^4
\AAA^a\BB^b$, with either, $a,b =0$ or $a =1, b=-1$, since the degree
of the zero and the pole divisor must be the same.

Case I: $a=b=0$\\
Since no prime which is a zero of $y$ ramifies in the extension
$\tilde K/ \tilde C_K(z)$, the divisor of $y$ in $\tilde C_K(z)$ is
also a fourth power of another divisor. In the rational function field
$\tilde C_K(z)$ every degree $0$ divisor is principal, so $y \in
(\tilde C_K(z))^4$.

Case II: $a=1, b= -1$\\
In this case, the divisor of $\frac{y}{z}$ is of the form $\EE^4$ and
hence $\frac{y}{z}= f^4$ for some $f \in \tilde C_K(z)$ by the same
argument as in Case I. Hence $y+z = u^4 + u$ can be rewritten as
$z \left( \frac{y}{z}+1 \right) = z(f+1)^4 = u^4 + u$. Since $f+1$ is a
rational function in $z$, we can rewrite this as
\begin{equation} \label{E:int}
z \left( \frac{f_1}{f_2} \right) ^4 = u^4 +u 
\end{equation}
where $ f_1, f_2$ are relatively prime polynomials in $\tilde
C_K[z]$, and $ f_2$ is monic. Equation (\ref{E:int}) shows: Any
valuation which is a pole of $u$ is either a pole of $z$ or a zero of
$f_2$. Let $\cc$ be a pole of $u$ which is a zero of $f_2$. Then,
since $f_2$ is a polynomial in $z$, $\cc$ is not a pole of $z$. So we
must have $|\ord_{\cc}f_2|= |\ord_{\cc} u|$. Hence $s:= f_2 \cdot u$
will have poles only at the valuations which are poles of $z$. Thus we
can rewrite (\ref{E:int}) in the form
\begin{equation}\label{F:int}
zf_1^4 + s^4 = s f_2^3.
\end{equation}
Furthermore, let $\cc$ be a zero of $f_2$. As pointed out above, $\cc$
is not a pole of $z$, so $\cc$ is not a pole of $s$. So we can deduce
that for a zero $\cc$ of $f_2$ we have $\ord_{\cc}(s^4 +
zf_1^4)=\ord_{\cc}(sf_2^3) \geq 3$.  Thus $\ord_{\cc}(d(s^4 +
zf_1^4))\geq 2$, so $\ord_{\cc}(f_1^4\,dz) \geq 2.$ Here $dz$ denotes
a K\"ahler differential. Since $\cc$ is unramified in the extension
$\tilde K/\tilde C_K(z)$, $\ord_{\cc}(dz)=0$. Hence $\ord_{\cc}(f_1^4)
\geq 2$, {\it i.e.}\ $f_1$ has a zero at $\cc$. Since $f_1$ and $f_2$
are relatively prime polynomials, this implies that $f_2$ has no
zeros, {\it i.e.}\ $f_2 = 1$. Hence $y$ is a polynomial in $z$.
Exactly the same argument applied to $\frac{1}{y}$ shows that
$\frac{1}{y}$ is a polynomial in $\frac{1}{z}$. Thus $y = z^{l}$ for
some $l \geq 0$ and $y+z = z^l +z = u^4 + u$. If $y=z$, we are done.
Otherwise this implies that all the poles of $y+z$ have order $l$ (the
poles of $z$ are simple), and also, that all the poles of $y+z$ are
divisible by 4.  Hence $4|l$.

So in both cases, Case I and Case II, we could deduce that either $y=z$
or that $y \in (\tilde C_K(z))^4$.
Since we assumed that $y \notin (\tilde C_K(z))^4$
this concludes the proof.
\end{proof}
%%%%%%%%%%%%%%%%%%%%%%%%%%%%%%%%%%%%%%%%%%%%%%%%%%%%%%%%%%%%%%%%%%
%%%%%%%%%%%%%%%%%%%%%%%%%%%%%%%%%%%%%%%%%%%%%%%%%%%%%%%%%%%%%%%%%
%%%%%%%%%%%%%%%%%%%%%%%%%%%%%%%%%%%%%%%%%%%%%%%%%%%%%%%%%%%%%%%%%

\begin{lemma}\label{tpower}

For all $c,c' \in V$ let $t_{c,c'}:= \frac{t+c}{t+c'}.$
  Let $w, v, u, u_{d,d'}, v_{d,d'}$ be elements of $K$
  such that $\forall c\in V\;\exists d \in V_c$ such that $\forall c'
  \in V \;\exists d' \in V_{c'}$
  such that the following equations are satisfied: 
\begin{gather}\label{F1:int}
  w +t = u^4 + u\\
  \label{F2:int} \frac{1}{w} + \frac{1}{t} = v^4 +v\\
  w_{d,d'}= \frac{w+d}{w+d'}\\
\label{F4:int} w_{d,d'}+t_{c,c'}= u^4_{d,d'} +u_{d,d'}\\
   \label{F3:int} \frac{1}{w_{d,d'}} +
   \frac{1}{t_{c,c'}}= v^4_{d,d'}+ v_{d,d'}
\end{gather}
Then $w = t^{4^s}$ for some natural number $s$.
\end{lemma}
%%%%%%%%%%%%%%%%%%%%%%%%%%%%%%%%%%%%%%%%%%%%%%%%%%%%%%%%%%%%%%%%%%%%
%%%%%%%%%%%%%%%%%%%%%%%%%%%%%%%%%%%%%%%%%%%%%%%%%%%%%%%%%%%%%%%%%%%%
%%%%%%%%%%%%%%%%%%%%%%%%%%%%%%%%%%%%%%%%%%%%%%%%%%%%%%%%%%%%%%%%%%%%
\begin{proof}
  Recall that the divisor of $t$ in $K$ is of the form $\pp/\qq$, and
  that $\PP$ and $\QQ$ are the primes of $C_K(t)$ lying below $\pp$
  and $\qq$, respectively. Thus the degree of $\QQ$ is one. Similarly,
  for all $c \in V$ the degree of the primes $\PP_c$ in $C_K(t)$ is
  one.  Hence $\QQ$ and all the $\PP_c$'s will remain prime in
  the constant field extension $\tilde C_K(t)/C_K(t)$. By Lemma~6.16
  in \cite{Sh2000} their factors will be unramified in the extension
  $\tilde K/\tilde C_K(t)$. Hence for all $c,c' \in V$, $t_{c,c'}$ has
  neither zeros nor poles at any prime ramifying in the extension
  $\tilde K/ \tilde C_K(t)$.

  In the second paragraph of the proof of Lemma~2.6 of \cite{Sh2000},
  pp.~471--472, translated to our notation, Shlapentokh proves that
  for some $c_0 \in V$ there exists a subset $V'$ of $V$ containing
  $n+1$ elements, not containing $c_0$, and such that
  for any $d_0 \in V_{c_0}$, for all $c' \in V'$, for any $d' \in
  V_{c'}$, $w_{d_0,d'}$ does not have zeros or poles at any prime
  ramifying in the extension $\tilde K/ \tilde C_K(t)$. Her argument
  uses the fact that there are exactly $r$ primes ramifying in the
  extension $\tilde K/\tilde C_K (t)$, and it does not use the
  characteristic of $K$, so the same proof works here. We have two
  cases to consider.
  
  Case I: $w \in C_K(t)$.\\
  If $w$ is in $C_K(t)$, then pick a $d_0 \in V_{c_0}$ and for some
  $c' \in V'$ pick a $d' \in V_{c'}$ such that (\ref{F4:int}) and
  (\ref{F3:int}) are satisfied. Then $w_{d_0,d'} \in C_K(t)$, and
  we can apply Lemma~\ref{basecase} to $t_{c_0,c'}$ instead of $t$,
  and to $w_{d_0,d'}$ to conclude that $w_{d_0,d'} =
  t_{c_0,c'} ^{4^s}$ for some $s \geq 0$. (Note that
  $C_K(t)=C_K(t_{c_0,c'})$.)  So
\[ \frac{w+{d_0}}{w+d'} = \left(
  t_{c_0,c'} \right ) ^{4^s}.\] If $s=0$, then we can check that
  $w=t$. (See the last part of Lemma~\ref{frac}.) Otherwise write
$1+ \frac{d_0+d'}{w+d'} = \left( t_{c_0,c'} \right )
  ^{4^s}.$
  
  Hence $w+d' = \left ( \frac{1}{t_{c_0,c'} + 1} \right )^{4^s}\cdot
  {(d_0+d')}$ . Since $(d_0+d')$ is an element of $C$ and hence a
  fourth power this implies that $w +d' \in (C_K(t))^4$, and hence $w
  \in (C_K(t))^4$, say $w = \tilde w^4$.  We can rewrite equations
  (\ref{F1:int}) and (\ref{F2:int}) as
\begin{gather}
  \tilde w +t = (u +\tilde w)^4 + (u + \tilde w)\\
  \frac{1}{\tilde w} + \frac{1}{t} = (v + \frac{1}{\tilde w})^4
  +(v+\frac{1}{\tilde w}).
\end{gather} 
Also
\[
w_{d,d'}=
\frac{w+d}{w+d'}=\frac{w+\tilde d^4}{w+ {{\tilde d}'}{}^4}=
\left( \frac{\tilde w + \tilde d}{\tilde w + \tilde d'}
\right)^4
\]
for $d\in V_c,d'\in V_{c'}$ and some suitable $\tilde d \in V_c$ and
$\tilde d' \in V_{c'}$. This lets us rewrite equations (\ref{F4:int})
and (\ref{F3:int}) in a similar fashion. So we can rewrite the
equations (\ref{F1:int}) through (\ref{F3:int}), and $\tilde w \in
C_K(t)$. Equation (\ref{F1:int}) and the fact that $t$ has only simple
zeros imply that $w \notin C_K$. Hence after finitely many iterations
we must be in the position where $s=0$.

Case II: $w \notin C_K(t).$\\
In this case we will derive a contradiction. $w \notin C_K(t)$ would
imply that $w_{d,d'} \notin C_K(t)$ for all $d$ and $d'$.

By putting $\alpha := u^2 +u$ we can rewrite equation (\ref{F1:int}) as
\begin{gather}
w+t = \alpha^2 + \alpha. \label{new1}
\end{gather}
Similarly by putting $\beta :=v^2 +v$,
 $\alpha_{d,d'}:=u_{d,d'}^2 + u_{d,d'}$,
 $\beta_{d,d'}:= v_{d,d'}^2 + v_{d,d'}$ we can
 rewrite (\ref{F2:int}),(\ref{F4:int}) and (\ref{F3:int}) as
\begin{gather}
  \frac{1}{w} + \frac{1}{t} = \beta^2 + \beta\\
  w_{d,d'}+t_{c,c'}=
  \alpha^2_{d,d'} +\alpha_{d,d'}\\
  \label{new2}\frac{1}{w_{d,d'}} +
  \frac{1}{t_{c,c'}}= \beta^2_{d,d'}+
  \beta_{d,d'}.
\end{gather} 

Let $c_0 \in V$ be as above. By the same argument as in \cite{Sh2000},
p.~472, with $p$ replaced by $2$, (\ref{new1}) through (\ref{new2})
imply that $\exists d_0 \in V_{c_0}$ such that $\forall c'\in
V'\;\exists d' \in V_{c'}$ such that the divisor of $w_{d_0,d'}$ is of
the form $\AAA^2\pp^a_{d'}\pp^b_{d_0}$. Here $\pp_{d_0}$ and
$\pp_{d'}$ are prime divisors, $a$ is either $-1$ or
$0$, and $b$ is either $1$ or $0$. Now the proof follows word for word
that of Lemma~2.6 in \cite{Sh2000}, p.~472 with $p$ replaced by $2$ to
prove that in this case $w = \tilde w^2$ with $\tilde w \in K$. For
this part of the proof we only used equations (\ref{new1}) through
(\ref{new2}).  Now we can rewrite equations (\ref{new1}) through
(\ref{new2}) with $w$ replaced by $\tilde w$.  Since $w \notin
C_K(t)$, $\tilde w \notin C_K(t)$. So we can keep replacing $w$ by its
square root over and over, contradicting that $w \in K^{2^s}$ only if
$s \leq [K:C_K(t)]$.  So $w = t^{4^s}$ for some $s \in \N$.
\end{proof}
%%%%%%%%%%%%%%%%%%%%%%%%%%%%%%%%%%%%%%%%%%%%%%%%%%%%%%%%%%%%%%%%%%
%%%%%%%%%%%%%%%%%%%%%%%%%%%%%%%%%%%%%%%%%%%%%%%%%%%%%%%%%%%%%%%%%%
%%%%%%%%%%%%%%%%%%%%%%%%%%%%%%%%%%%%%%%%%%%%%%%%%%%%%%%%%%%%%%%%%%
%%%%%%%%%%%%%%%%%%%%%%%%%%%%%%%%%%%%%%%%%%%%%%%%%%%%%%%%%%%%%%%%%%
\begin{cor}\label{diophantine0}
The set $S_1:=\{\, t^{4^s}:s \in \N\,\}$ is
  diophantine over $K$.
\end{cor}
\begin{proof}
  Lemma~\ref{tpower} shows that an element $w \in K$ satisfying the
  equations (\ref{F1:int}) through (\ref{F3:int}) must be of the form
  $w = t^{4^s}$ for some $s \in \N$. What we have left to show is that
  if $w = t ^{4^s}$ for some $s \in \N$, then we can satisfy equations
  (\ref{F1:int}) through (\ref{F3:int}).  If $w=t$, let $u=0$ and
  $v=0$. For the general case we use the fact that for any $x \in K$
  and any $s \in \N$ we have
\[
x ^{4^s}-x = (x^{4^{s-1}}+ x^{4^{s-2}}+ \dots + x)^4-(x^{4^{s-1}}+
x^{4^{s-2}}+ \dots + x).
\] 
So if $w = t^{4^s}$ with $s \geq 1$, let $u = t^{4^{s-1}} + \dots + t^4
+t$. For $v$ take
\[
\frac{1}{t^{4^{s-1}}}+ \dots + \frac{1}{t^4}+\frac{1}{t}.
\]
Now fix $c\in V$.  To satisfy the other equations we can use the same
argument, if we can show that $\exists d \in V_c$ such that $\forall
c' \in V \exists d' \in V_{c'}$ such that $w_{d,d'} =
(t_{c,c'})^{4^s}$. This is done in Corollary~2.7 in \cite{Sh2000}.
\end{proof}
\begin{cor}\label{diophantine1}
  The set $S = \{\,t^{2^s}: s \in \N \,\}$ is diophantine over $K$.
\end{cor}
\begin{proof}
This follows from the fact that 
\[
w \in S \iff  (w \in S_1 {\mbox{ or }} \exists z
\in K \,(z^2 =w {\mbox { and }} z\in S_1)). 
\] 
\end{proof}
%%%%%%%%%%%%%%%%%%%%%%%%%%%%%%%%%%%%%%%%%%%%%%%%%%%%%%%%%%%%%%%%%%%
%%%%%%%%%%%%%%%%%%%%%%%%%%%%%%%%%%%%%%%%%%%%%%%%%%%%%%%%%%%%%%%%%%%
%%%%%%%%%%%%%%%%%%%%%%%%%%%%%%%%%%%%%%%%%%%%%%%%%%%%%%%%%%%%%%%%%%%

%%%%%%%%%%%%%%%%%%%%%%%%%%%%%%%%%%%%%%%%%%%%%%%%%%%%%%%%%%%%%%%%%%
%%%%%%%%%%%%%%%%%%%%%%%%%%%%%%%%%%%%%%%%%%%%%%%%%%%%%%%%%%%%%%%%%%
\subsection{\large {The set $T=\left\{(x,w) \in K^2: \exists s \in 
      \N, w = {\left ( \frac{x^2 + t^2+t}{x^2 +t} \right
        )}^{2^s}\right\}$ is diophantine over $K$.}}
\begin{lemma}\label{simplezero}
  Let $x \in K$. Let $t$ be as above, {\it i.e.}\ $\tilde K / \tilde C
  _K(t)$ is separable and the divisor of $t$ is of the form
  $\frac{\pp}{\qq}$. Let $u = \frac{x^2 + t^2 +t}{x^2 +t}$, and let $a
  \in C$, $a \neq 1$. Then $u+a$ has only simple zeros and simple
  poles, except possibly for zeros at $\pp,\qq$ or primes ramifying in
  the extension $\tilde K/\tilde C_K(t)$.
\end{lemma}
\begin{proof}
  First we will show that the zeros of $u+a$ away from the ramified
  primes and $\pp$ and $\qq$ are simple: By Lemma~4.4 in \cite{Sh96}
  it is enough to show that $u+a$ and
  $\frac{du}{dt}$ do not have common zeros. We have
\begin{eqnarray*}
\frac{d(u+a)}{dt}&=&\frac{(x^2 + t^2 +t)+(x^2
 +t)}{(x^2+t)^2}=\frac{t^2}{(x^2+t)^2}\;\mbox{ and }\\
u+a &=& \frac{x^2 + t^2 +t}{x^2+t}+a = 1+a + \frac{t^2}{x^2 +t}.
\end{eqnarray*}
Suppose $\cc$ is a zero of $d(u+a)/dt$ satisfying the above
conditions. Then $\cc$ is not a zero of $t$, so $\cc$ must be a pole
of $x^2 +t$, {\it i.e.}\ a pole of $x$. If $\cc$ is a pole of $x$,
then it is a zero of $\frac{t^2}{x^2 + t}$, and hence not a zero of $1
+ a + \frac{t^2}{x^2 +t}$. Hence $d(u+a)/dt$ and $u+a$ have no zeros
in common, except possibly the ones mentioned above.

We will now show that all poles at above described valuations are
simple: Since $u$ and $u+a$ have the same poles, it is enough to show
that the poles of $u$ are simple. $u$ has simple poles if and only if
the zeros of $u^{-1}$ are simple. So we'll show that the zeros of
$v=u^{-1}$ are simple by doing exactly the same thing as above. Let
$v:= u^{-1} = \frac{x^2 +t}{x^2 + t^2 +t}$. Then
\begin{eqnarray*}
\frac{dv}{dt}&=& \frac{(x^2 +t +x^2 +t^2 +t)}{(x^2 +t^2 +t)^2}=
\frac{t^2}{(x^2 +t^2 +t)^2}\;\mbox{ and }\\
v &=& \frac{x^2 +t}{x^2 + t^2 +t}=1 - \frac{t^2}{(x^2+t^2+t)^2}.
\end{eqnarray*}
Again let $\cc$ be a zero of $dv/dt$ satisfying the above conditions.
Again $\cc$ has to be a pole of $x$. So $\cc$ is a zero of $v$, but
not a zero of $1 - \frac{t^2}{x^2 +t^2 +t}$, since $\cc$ is not a zero
or pole of $t$. Hence all the zeros of $u^{-1}$ are simple except
possibly for the ones mentioned above.
\end{proof}

%%%%%%%%%%%%%%%%%%%%%%%%%%%%%%%%%%%%%%%%%%%%%%%%%%%%%%%%%%%%%%%%%
%%%%%%%%%%%%%%%%%%%%%%%%%%%%%%%%%%%%%%%%%%%%%%%%%%%%%%%%%%%%%%%%%
%%%%%%%%%%%%%%%%%%%%%%%%%%%%%%%%%%%%%%%%%%%%%%%%%%%%%%%%%%%%%%%%%
%%%%%%%%%%%%%%%%%%%%%%%%%%%%%%%%%%%%%%%%%%%%%%%%%%%%%%%%%%%%%%%%%%
%%%%%%%%%%%%%%%%%%%%%%%%%%%%%%%%%%%%%%%%%%%%%%%%%%%%%%%%%%%%%%%%%%
%%%%%%%%%%%%%%%%%%%%%%%%%%%%%%%%%%%%%%%%%%%%%%%%%%%%%%%%%%%%%%%%%%

\begin{lemma}\label{upower}
  Let $x,v \in K^*$, let $u := \frac{x^2 + t^2 +t}{x^2 + t}$.  For
  all $c,c' \in V$, $g\in\{-1,1\}$ let
\[
  u_{c,c',g} := \frac{u^g + c}{u^g + c'}.\]
For $d \in V_c, d'\in V_{c'}, g \in\{-1,1\}$ let
\[  v_{d,d',g} := \frac{v^g + d}{v^g + d'}.
\]
In addition assume that $\forall c\in V \;\exists d \in V_c$ such that
$\forall c' \in V \;\exists d' \in V_{c'}$ such that the following
equations hold for for $e,g \in \{-1,1\}$, and some $s \in \N$:
\begin{gather}\label{G1:int}
 v^{e}_{d,d',g}+ u^{e}_{c,c',g}=
 \sigma^4_{d,d',e,g}+\sigma_{d,d',e,g}\\
 \label{G3:int} v^{2e}_{d,d',g}t^{4^s}+
 u^{2e}_{c,c',g}t = \lambda_{d,d',e,g}^4 +
 \lambda_{d,d',e,g}\\
 \label{G5:int} (u^g + c)^e + (v^g + d)^e =
 \mu^4_{d,e,g}+\mu_{d,e,g}
\end{gather}
Then for some natural number $m$, $v=u^{4^m}$.
\end{lemma}
\begin{proof}
It is sufficient to prove that the result is valid in $\tilde K :=
\tilde C_K K$, so we work in $\tilde K$.
We will first prove the following

{\it {Claim}}: For all $c,c' \in V,g\in \{-1,1\}, \;
u_{c,c',g}$ has no multiple zeros or poles except possibly at the
primes ramifying in $\tilde K/\tilde C_K(t)$ or $\pp$ or $\qq$.

{\it Proof of Claim: }By Lemma~\ref{frac} we have that for all
$c,c',g$ as above all the poles of $u_{c,c',g}$ are zeros of $u^g +
c'$ and all the zeros of $u_{c,c',g}$ are zeros of $u^g + c$. By
Lemma~\ref{simplezero} and by assumption on $c$ and $c'$, all the
zeros of $u^g +c'$ and $u^g + c$ are simple, except possibly for zeros
at $\pp$, $\qq$ or primes ramifying in the extension $\tilde K/ \tilde
C_K(t)$. This proves the claim.

Also since $\frac{du}{dt} \neq 0$, $u$ is not a second power in
$\tilde K$.
We will show the following: (a) If $s=0$, then $u=v$, and (b) if $s
>0$, then $v$ is a fourth power of some element in $\tilde K$.

Case~(a): Suppose that $s =0$, and set $g =1$.

Again, using Shlapentokh's argument in Lemma~2.6 of \cite{Sh2000}
there exists $c_0 \in V$ and $V'\subseteq (V-\{c_0\})$ containing
$n+1$ elements, such that for all $d_0 \in V_{c_0}$, for all $c' \in
V'$, and for all $d' \in V_{c'}$, $u_{c_0, c',1}$ and $v_{d_0,d',1}$ have
no zeros or poles at the primes of $\tilde K$ ramifying in the
extension $\tilde K/ \tilde C(t)$ or at $\pp$ or $\qq$. For indices
selected in this way, all the poles and zeros of $u_{c_0,c',1}$ are
simple. Pick $d_0 \in V_{c_0}$
and for all $c' \in V'$ pick $d'\in V_{c'}$ such that equations
(\ref{G1:int}) and (\ref{G3:int}) are satisfied.  Equations
(\ref{G3:int}) and (\ref{G1:int}) imply:
\begin{gather}\label{H1:int}
  v^2_{d_0,d',1}t^{4^s} + u^2_{c_0,c',1}t =
  \lambda^4_{d_0,d',1,1}+ \lambda_{d_0,d',1,1}\\
\label{H3:int}
   v^2_{d_0,d',1}+u^2_{c_0,c',1} =
  \sigma^8_{d_0,d',1,1}+ \sigma^2_{d_0,d',1,1}.
\end{gather}

From (\ref{H1:int}) and (\ref{H3:int}) we obtain (since $s=0$)
\[
\lambda^4_{d_0,d',1,1}+ \lambda_{d_0,d',1,1} =
t( \sigma^8_{d_0,d',1,1}+ \sigma^2_{d_0,d',1,1}).
\]
All the poles of $\lambda_{d_0,d',1,1}$ and
$\sigma_{d_0,d',1,1}$ are poles of $u_{c_0,c',1}$,
$v_{d_0,d',1}$ or $t$, and thus are not at any valuation
ramifying in the extension $\tilde K/\tilde C_K(t)$. By
Lemma~\ref{sigma} applied to $\sigma =
\sigma^2_{d_0,d',1,1}$ and (\ref{H3:int})
\[
v^2_{d_0,d',1} + u^2_{c_0,c',1} =0.
\]
Thus $v_{d_0,d',1} = u _{c_0,c',1}$. From here on the proof
is word for word like the proof of Lemma~2.10 in \cite{Sh2000}, top of
page 477, showing that if $s=0$, then $u=v$.

Case (b): Suppose now that $s >0$. Again set $g=1$. Let $c_0$ and $V'$
be selected as above.

Again we can pick $d_0 \in V_{c_0}$ and for all $c' \in V'$ we can
pick $d' \in V_{c'}$ such that equations (\ref{G1:int}) through
(\ref{G5:int}) are sa\-tis\-fied and such that the corresponding
$u_{c_0,c',1}$ and $v_{d_0,d',1}$ do not have zeros or poles at the
primes of $\tilde K$ ramifying in the extension $\tilde K/ \tilde
C(t)$ or at $\pp$ or $\qq$.  We can use the same argument as in
Lemma~\ref{basecase} to show that either (i) for all $d'$ chosen as
above the divisor of $v_{d_0,d',1}$ in $\tilde K$ is a fourth power of
another divisor or (ii) for some $c' \in V'$ and some $d' \in V_{c'}$
and some prime $\ttt$ not ramifying in $\tilde K/\tilde C(t)$ and not
equal to $\pp$ or $\qq$, $\ord_{\ttt}v_{d_0,d'} \in \{1,-1\}$. In
Case~(i), because of our choice of the $v_{d_0,d'}$'s and Proposition
\ref{Lemma2.9}, a short calculation shows that $v \in \tilde K^4$:
\begin{eqnarray*}
v_{d_0,d',1}^{-1}&=& \frac{v+d_0}{v+d'}=1+
\frac{(d'+d_0)}{v+d_0}\\
 &=&(d'+d_0)
(\frac{1}{d'+d_0}+\frac{1}{v+d_0}),
\end{eqnarray*}
where $d_0 \in V_{c_0}$ is fixed, and we have $d' \in V_{c'}$ , and all
$d'$ are distinct. Also $V'$ contains $n+1$ elements, so by
Proposition \ref{Lemma2.9} applied to $\frac{1}{v+d'}$ we have that
$\frac{1}{v+d'} \in \tilde K^{4}$. This implies that $v \in
\tilde K^{4}$.
This finishes Case (i).

So assume now that we are in Case (ii): Without loss of generality,
assume that $\ttt$ is a pole of $v_{d_0,d',1}$ (and hence
neither a zero nor a pole of $t$).

Again look at equations (\ref{H1:int}) and (\ref{H3:int}).
Since $t$ does not have a pole or a zero at $\ttt$ and since the right
hand sides of equations (\ref{H1:int}) and (\ref{H3:int}) only have
poles of order $\geq 4$,
\begin{gather}
  \ord_{\ttt}(v^2_{d_0,d',1}t^{4^s} + u^2_{c_0,c',1}t) =
  \ord_{\ttt}(\lambda^4_{d_0,d',1,1}+
  \lambda_{d_0,d',1,1})\geq 0 \mbox{ and }\\
  \ord_{\ttt}(v^2_{d_0,d',1}+u^2_{c_0,c',1}) =\ord_{\ttt}
  (\sigma^8_{d_0,d',1,1}+
  \sigma^2_{d_0,d',1,1})\geq 0
\end{gather}
Thus 
\begin{eqnarray*}
\ord_{\ttt}v^2_{d_0,d',1}(t^{4^s}+t)\geq 0
\end{eqnarray*}
Hence it follows that $\ord_{\ttt}(t^{4^s} + t ) \geq 2
|\ord_{\ttt}v_{d_0,d',1}| $. But in $\tilde C_K(t)$ all the zeros
of $t^{4^s}+t$ are simple. So this function can have multiple zeros
only at primes ramifying in the extension $\tilde K/\tilde C_K(t)$.
But by assumption $\ttt$ is not one of these primes, and so we have a
contradiction unless $v \in \tilde K^4$.
This shows that if $s>0$, then equations (\ref{G1:int}) through
(\ref{G5:int}) can be rewritten in the same fashion as in Lemma
\ref{tpower} with $v$ replaced by its fourth roots, and in
(\ref{G3:int}) $s$ is replaced by $s-1$. Therefore, after finitely
many iterations of this rewriting procedure we will be in the case of
$s=0$, which was treated in Case~(a). Hence, for some natural number
$m$, $v=u^{4^m}$.
\end{proof}
%%%%%%%%%%%%%%%%%%%%%%%%%%%%%%%%%%%%%%%%%%%%%%%%%%%%%%%%%%%%%%%%%%%%
%%%%%%%%%%%%%%%%%%%%%%%%%%%%%%%%%%%%%%%%%%%%%%%%%%%%%%%%%%%%%%%%%%%
%%%%%%%%%%%%%%%%%%%%%%%%%%%%%%%%%%%%%%%%%%%%%%%%%%%%%%%%%%%%%%%%%%%
%%%%%%%%%%%%%%%%%%%%%%%%%%%%%%%%%%%%%%%%%%%%%%%%%%%%%%%%%%%%%%%%%%%%
%%%%%%%%%%%%%%%%%%%%%%%%%%%%%%%%%%%%%%%%%%%%%%%%%%%%%%%%%%%%%%%%%%%
\begin{cor}
  The set $T_1:=\left\{(x,w) \in K^2: \exists s \in \N, w = {\left (
      \frac{x^2 + t^2+t}{x^2 +t} \right )}^{4^s}\right\}$ is diophantine
  over $K$.
\end{cor}
\begin{proof}
 Let $x \in K$, and let $u = \frac{x^2 + t^2 +t}{x^2 +t}$.
  Lemma~\ref{upower} shows that an element $v \in K$ satisfying
  equations (\ref{G1:int}) through (\ref{G5:int}) must be of the form
  $v = u^{4^k}$ for some $k \in \N$. So we have to show now that if $v
  = u^{4^k}$ for some $k \in \N$, then equations (\ref{G1:int})
  through (\ref{G5:int}) can be satisfied.  The proof of this is
  almost identical to Corollary~\ref{diophantine0}.
\end{proof}
\begin{cor}\label{diophantine2}
  The set $T:=\left\{(x,w) \in K^2: \exists s \in \N, w = {\left (
      \frac{x^2 + t^2+t}{x^2 +t} \right )}^{2^s}\right\}$ is diophantine
  over $K$.
\end{cor}
\begin{proof}
This follows from the fact that 
\[
(x,w) \in T \iff (x,w) \in T_1 \mbox{ or } \exists z \in K
\,(z^2 = w \mbox{ and } (x,z) \in T_1).
\]
\end{proof}
%%%%%%%%%%%%%%%%%%%%%%%%%%%%%%%%%%%%%%%%%%%%%%%%%%%%%%%%%%%%%%%%%%%%%%
%%%%%%%%%%%%%%%%%%%%%%%%%%%%%%%%%%%%%%%%%%%%%%%%%%%%%%%%%%%%%%%%%%%
%%%%%%%%%%%%%%%%%%%%%%%%%%%%%%%%%%%%%%%%%%%%%%%%%%%%%%%%%%%%%%%%%%
%%%%%%%%%%%%%%%%%%%%%%%%%%%%%%%%%%%%%%%%%%%%%%%%%%%%%%%%%%%%%%%%%%
\begin{theorem}\label{THEEND}
  The set $\{(x,y) \in K^2 : \exists s \in \N, y= x ^{2^s}\}$ is
  diophantine over $K$.
\end{theorem}
\begin{proof}
  By Corollary~\ref{diophantine1}, Corollary~\ref{diophantine2} and
  the remark after Theorem~\ref{mainpart}, the sets $S,S',T$, and $T'$
  are diophantine. Together with Theorem~\ref{mainpart} this finishes
  the proof.
\end{proof}
\section{Appendix}
%%%%%%%%%%%%%%%%%%%%%%%%%%%%%%%%%%%%%%%%%%%%%%%%%%%%%%%%%%%%%%%%%%%%
%%%%%%%%%%%%%%%%%%%%%%%%%%%%%%%%%%%%%%%%%%%%%%%%%%%%%%%%%%%%%%%%%%%
%%%%%%%%%%%%%%%%%%%%%%%%%%%%%%%%%%%%%%%%%%%%%%%%%%%%%%%%%%%%%%%%%%
%%%%%%%%%%%%%%%%%%%%%%%%%%%%%%%%%%%%%%%%%%%%%%%%%%%%%%%%%%%%%%%%%%%
In the appendix we give proofs for Proposition~\ref{Lemma2.9} and
Lemma~\ref{sigma}. Both were used in Lemma~\ref{upower}.
\begin{lemma}
  Let $F/G$ be a finite extension of fields of positive
  characteristic $p$. Let $\alpha \in F$ be such that all the
  coefficients of its monic irreducible polynomial over $G$ are
  in $G^{p^2}$. Then $\alpha \in F^{p^2}$.
\end{lemma}
\begin{proof}
  This is Lemma~2.1 in \cite{Sh2000} with $p$ replaced by $p^2$, and
  the same proof works here.
\end{proof}
\begin{cor}\label{help}
  Let $F/G$ be a finite separable extension of fields of positive
  characteristic $p$. Let $[F:G]=n$. Let $x \in F$ be such that
  $F=G(x)$, and such that for distinct $a_0,\dots,a_n \in G$,
  $\Norm_{F/G}(a_i^{p^2}-x)=y_i^{p^2}$ with $y_i \in G$. Then $x \in
  F^{p^2}$.
\end{cor}
\begin{proof}
  This is very similar to Lemma~2.2 in \cite{Sh2000}. Let $H(T) = A_0
  + A_1T + \dots + T^n$ be the irreducible polynomial of $x$ over $G$.
  Then $H(a_i^{p^2})=y_i^{p^2}$ for $i \in \{0,\dots,n\}$. This gives
  us the following linear system of equations:
\[
\left( \begin{array}{ccccc}
    1 & a_0^{p^2} & \dots & a_0^{p^2 (n-1)} &a_0^{p^2n}\\
    \vdots & \vdots & \ddots & \vdots & \vdots \\
    1& a_n^{p^2} & \dots & a_n^{p^2(n-1)}&a_n^{p^2n}
\end{array}\right)
\left( \begin{array}{c}
A_0\\ \vdots \\ 1
\end{array} \right) = 
\left(
\begin{array}{c}
y_0^{p^2}\\
\vdots\\
y_n^{p^2}
\end{array}
\right)
\]
We can use Cramer's rule to solve the system and to conclude that
$A_i \in G^{p^2}$ for all $i$. Then by the previous lemma, $x \in
F^{p^2}$.

\end{proof}
Now we can apply the corollary to our situation:
\begin{prop}\label{Lemma2.9}
  Let $v \in \tilde K$, and assume that for some distinct $a_0, \dots,
  a_n \in C$, the divisor of $v+a_i$ is of the form
  $\mathcal{D}_i^{p^2}$ for divisors $\mathcal{D}_i$ of $\tilde K$,
  $i=0,\dots,n$. Moreover, assume that for all $i$, $v+a_i$ does not
  have zeros or poles at any prime ramifying in the extension $\tilde
  K/ \tilde C_K(t)$. Then $v \in \tilde K^{p^2}$.
\end{prop}
\begin{proof}
  This is almost the same as Lemma~2.9 in \cite{Sh2000}: First assume
  that $v \in \tilde C_K(t)$. Since $v+a_i$ does not have any zeros or
  poles at primes ramifying in the extension $\tilde K/ \tilde
  C_K(t)$, the divisor of $v+a_i$ in $\tilde C_K(t)$ is of the form
  $\mathcal{E}_i^{p^2}$. In $\tilde C_K(t)$ every divisor of degree
  zero is principal, so $v+a_i \in (\tilde C_K(t))^{p^2}$ and hence $v
  \in (\tilde C_K(t))^{p^2}$. Therefore $v \in \tilde K^{p^2}$.

  So now assume that $v \notin \tilde C_K(t)$. From our assumption on
  $v+a_i$ it follows that in $\tilde C_K(t,v)$ the divisor of $v+a_i$
  is a $p^2$ power of another divisor. Since the divisor of
  $\Norm_{\tilde C_K(t,v)/ \tilde C_K(t)}(v+a_i)$ is equal to the
  corresponding norm of the divisor of $(v+a_i)$, it follows that the
  divisor of the $\tilde C_K(t,v)/\tilde C_K(t)$ norm of $(v+a_i)$ is
  of the form $\mathcal{N}_i^{p^2}$, and hence $\Norm_{\tilde
    C_K(t,v)/\tilde C_K(t)}(v+a_i) \in (\tilde C_K(t))^{p^2}$. Now apply
  Corollary~\ref{help} with $G=\tilde C_K(t)$ and $F=\tilde C_K(t,v)$.
\end{proof}
\begin{lemma}\label{sigma}
  Let $\sigma, \mu \in K$. Assume that all the primes that are poles
  of $\sigma$ or $\mu$ do not ramify in the extension $\tilde K/\tilde
  C_K(t)$. Moreover assume that
\begin{equation}\label{L:int}
t (\sigma ^4 + \sigma)= \mu^4 + \mu.
\end{equation}
Then $\sigma^4 + \sigma = \mu^4 + \mu =0.$
\end{lemma}
\begin{proof}
  Let $\AAA$, $\BB$ be effective divisors of $K$, relatively prime to
  each other and to $\pp$ and $\qq$, such that the divisor of $\sigma$
  is of the form $\frac{\AAA}{\BB}\pp^i\qq^k$, where $i$ and $k$ are
  integers.
  
  {\it Claim 1:} For some effective divisor $\CC$ relatively prime to
  $\BB$, $\pp$ and $\qq$, some integers $j,m$, the divisor of $\mu$ is
  of the form $\frac{\CC}{\BB}\pp^j\qq^m$.

{\it Proof of Claim 1:} Let $\ttt$ be a pole of $\mu$ such that $\ttt
\neq \pp$ and $\ttt \neq \qq$. Then
\[
0 > 4\, \ord_{\ttt}\mu = \ord_{\ttt}(\mu ^4 + \mu) =
\ord_{\ttt}(t(\sigma ^4 + \sigma)) = \ord_{\ttt} (\sigma ^4 +
\sigma) = 4 \,\ord_{\ttt}\sigma.
\]
Conversely, let $\ttt$ be a pole of $\sigma$ such that $\ttt \neq
\pp$ and $\ttt \neq \qq$. Then 
\[
0 > 4 \,\ord_{\ttt}\sigma = \ord_{\ttt}(\sigma^4 + \sigma) =
\ord_{\ttt}(t(\sigma^4 + \sigma)) = \ord_{\ttt}(\mu ^4 + \mu) = 4
\,\ord_{\ttt}\mu.
\]
This proves the claim.

By the Strong Approximation Theorem there exists $b \in K^*$ such
that the divisor of $b$ is of the form $\frac{\BB\DD}{\qq^l}$, where
$\DD$ is an effective divisor relatively prime to $\AAA$, $\CC$, $\pp$
and $\qq$ and $l$ is a natural number.

{\it Claim 2:}
\begin{gather*}
  b \sigma = s_1 t^i\\
  b \mu = s_2 t^j ,
\end{gather*}
where $ s_1, s_2$ are integral over $ C_K[t]$ and have zero divisors
relatively prime to $ \pp$ and $\BB$.

{\it Proof of Claim 2:} The divisor of $b \sigma$ is
\[
\frac{\BB \DD}{\qq^l}\frac{\AAA}{\BB}\pp^i\qq^k = \DD \AAA \pp^i
\qq^{k-l} = (\DD \AAA \qq^{k-l+i})\left( \frac{\pp^i}{\qq^i} \right).
\]
Thus the divisor of $s_1 := b \sigma / t^i$ is of the form
$\DD\AAA\qq^{k-l+i}$. Therefore $\qq$ is the only pole of $s_1$, so
$s_1$ is integral over $C_K[t]$. By construction $\AAA$ and $\DD$ are
relatively prime to $\pp$ and $\BB$. A similar argument applies to
$s_2:=b\mu/t^j.$ This proves the claim.

Multiplying (\ref{L:int}) by $b^4$ we obtain the following equation
(using the claim):
\begin{gather}\label{M:int}
t(s_1^4t^{4i} +b^3s_1t^i)= s_2^4t^{4j}+ b^3s_2t^j.
\end{gather}
Suppose $i <0$. Then the left side of (\ref{M:int}) has a pole of
order $|4i +1|$ at $\pp$. This would imply that $j <0$, and the right
side has a pole of order $|4j|$ at $\pp$, contradiction. Thus we can
assume that $i,j$ are both nonnegative. We can rewrite
(\ref{M:int}) as 
\[
(s_1^4 t^{4i+1} + s_2 ^4 t^{4j}) = b^3(s_1 t^{i+1}+s_2 t^j).
\]
Let $\ttt$ be any prime factor of $\BB$ in $\tilde K$. Then $\ttt$
does not ramify in the extension $\tilde K / C_{\tilde K}(t)$ by our
assumption on $\sigma$. Also $\ttt$ is not a pole of $s_1, s_2$ or $t$.

Thus $$\ord_{\ttt}(s_1^4 t^{4i+1}+s_2^4t^{4j}) =
\ord_{\ttt}(b^3(s_1t^{i+1}+ s_2t^j)) \geq 3.$$ 
We have 
\[
0 <\ord_{\ttt}(d(s_1^4t^{4i+1}+s_2^4t^{4j}))=\ord_{\ttt}(s_1^4
\,d(t^{4i+1}))= \ord_{\ttt}(s_1^4)+\ord_{\ttt}(d(t^{4i+1})).
\]
Since $\ttt$ is unramified in the extension $\tilde K/\tilde C_K(t)$
and since $\ttt$ is not a zero or a pole of $t$,
$\ord_{\ttt}(d(t^{4i+1}))=0.$
So $s_1$ has a zero at $\ttt$. This, however, is impossible, because
$\ttt$ is a prime factor of $\BB$, but the zero divisor of $s_1$ is
relatively prime to $\BB$.  So $\BB$ must be the trivial divisor. This
implies that in (\ref{L:int}) all the functions are integral over
$C_K[t]$, {\it i.e.}\ they can have poles at $\qq$ only. So if $\mu$
is not constant, it must have a pole at $\qq$. But then the left side
of (\ref{L:int}) has a pole at $\qq$ of odd order, while the right
side of (\ref{L:int}) has a pole at $\qq$ of even order, which is a
contradiction.

Thus $\mu$ must be a constant. 
But if a function $h \in K$ is integral over $C_K[t]$, and $t \cdot h$ is
constant, then $h=0$.
Thus $\sigma^4 + \sigma = 0$. Then $\mu^4 + \mu =0$ also.
\end{proof}
%%%%%%%%%%%%%%%%%%%%%%%%%%%%%%%%%%%%%%%%%%%%%%%%%%%%%%%%%%%%%%%%%%
%%%%%%%%%%%%%%%%%%%%%%%%%%%%%%%%%%%%%%%%%%%%%%%%%%%%%%%%%%%%%%%%%%%%%
%%%%%%%%%%%%%%%%%%%%%%%%%%%%%%%%%%%%%%%%%%%%%%%%%%%%%%%%%%%%%%%%%%%%%
%%%%%%%%%%%%%%%%%%%%%%%%%%%%%%%%%%%%%%%%%%%%%%%%%%%%
%%%%%%%%%%%%%%%%%%%%%%%%%%%%%%%%%%%%%%%%%%%%%%%%%%%%

\end{document}